# ON THE CONVERGENCE TO EQUILIBRIUM OF KAC'S RANDOM WALK ON MATRICES

By Roberto Imbuzeiro Oliveira[1]

*Instituto de Matemática Pura e Aplicada (IMPA)*

We consider Kac's random walk on $n$-dimensional rotation matrices, where each step is a random rotation in the plane generated by two randomly picked coordinates. We show that this process converges to the Haar measure on $SO(n)$ in the $L^2$ transportation cost (Wasserstein) metric in $O(n^2 \ln n)$ steps. We also prove that our bound is at most a $O(\ln n)$ factor away from optimal. Previous bounds, due to Diaconis/Saloff-Coste and Pak/Sidenko, had extra powers of $n$ and held only for $L^1$ transportation cost.

Our proof method includes a general result of independent interest, akin to the path coupling method of Bubley and Dyer. Suppose that $P$ is a Markov chain on a Polish length space $(M,d)$ and that for all $x,y \in M$ with $d(x,y) \ll 1$ there is a coupling $(X,Y)$ of one step of $P$ from $x$ and $y$ (resp.) that contracts distances by a $(\xi + o(1))$ factor on average. Then the map $\mu \mapsto \mu P$ is $\xi$-contracting in the transportation cost metric.

**1. Introduction.** Around 50 years ago Kac [7] introduced a one-dimensional toy model of a Boltzmann gas. It is a discrete-time Markov process whose state at a time $t \in \{0,1,2,3,\ldots\}$ is a vector

$$v(t) = (v_1(t), \ldots, v_n(t)) \in \mathbb{R}^n,$$

corresponding to the velocities of $n$ interacting particles of equal mass. At each time $t$, a uniformly distributed pair $1 \leq i_t < j_t \leq n$ and a uniform angle $\theta_t \in [0, 2\pi]$ are chosen independently. This choice corresponds to a collision between particles $i_t, j_t$ whose velocities are changed to new values

$$v_{i_t}(t+1) = \cos\theta_t v_{i_t}(t) + \sin\theta_t v_{j_t}(t),$$
$$v_{j_t}(t+1) = -\sin\theta_t v_{i_t}(t) + \cos\theta_t v_{j_t}(t),$$

Received June 2007; revised June 2008.
[1]Supported in part by a *Bolsa de Produtividade em Pesquisa* from Conselho de Desenvolvimento Técnico e Científico (CNPq), Brazil.

whereas the other velocities are kept the same. This prescription for the new velocities implies that the *total kinetic energy*

$$E(t) \equiv \sum_{k=1}^{n} v_k(t)^2$$

is conserved.

For each time step $t$,

$$v(t+1) = R(i_t, j_t, \theta_t)v(t),$$

where $R(i_t, j_t, \theta_t)$ is a rotation by $\theta_t$ of the plane generated by the coordinates $i_t$ and $j_t$ in $n$-dimensional space. Two related processes have been studied in the literature under the heading of "Kac's random walk":

- Suppose $E(0) = 1$. Then the evolution of $v(0), v(1), v(2), v(3), \ldots$ corresponds to an ergodic Markov chain over the $(n-1)$-dimensional sphere $S^{n-1} \subset \mathbb{R}^n$, with uniform invariant distribution. This is the model originally considered by Kac [7] in his investigations of the foundations of Statistical Mechanics. See [4, 6, 16] and the references therein for more works in similar directions.
- One might also consider the random walk on matrices determined by choosing some $X(0) \in SO(n)$, the set of $n \times n$ rotation matrices, and then setting

$$X(t+1) = R(i_t, j_t, \theta_t)X(t), \qquad t \geq 0.$$

This is a discrete-time ergodic random walk on $SO(n)$ whose stationary distribution is a Haar measure on $SO(n)$. This process has also been extensively studied, both for its intrinsic interest and as a sampling algorithm—indeed, a "Gibbs sampler" [5]—for a Haar measure. Interestingly, this process is featured in Hastings' seminal 1970 paper on Markov chain Monte Carlo [8]. See [1, 4, 5, 12] for more details.

The question arises of how fast Kac's random walk on $SO(n)$ converges to equilibrium. This question may be posed in different forms. Convergence of density functions to equilibrium is very well understood: Janvresse [6] obtained the first bound of optimal magnitude $\Theta(n^{-1})$ on the $L^2$ spectral gap of the chain on $S^{n-1}$. Carlen, Carvalho and Loss [4] obtained the exact spectral gap for both processes. Finally, Maslen [10] computed the entire spectrum for both processes.

Convergence to equilibrium in total variation also occurs, as shown by Diaconis and Saloff-Coste [5] who obtained a very poor $e^{O(n^2)}$ mixing time bound for convergence in total variation of the matrix process. We cannot improve on this bound, but note that total variation is perhaps too stringent a notion of convergence for simulations (as it is sensitive to errors at arbitrarily small scales), whereas convergence of densities is too weak (e.g., when one starts from a discrete distribution).

We consider an intermediate notion of convergence to equilibrium based on *transportation cost*. Given a metric space $(M,d)$ and two probability measures $\mu$, $\nu$ over the Borel $\sigma$-field of $M$, the $L^p$ *transportation cost (or Wasserstein) distance between $\mu$ and $\nu$* is

$$W_{d,p}(\mu,\nu) = \inf\{(\mathbb{E}[d(X,Y)^p])^{1/p} : (X,Y) \text{ is a pair}$$
$$\text{of random variables coupling } (\mu,\nu)\}$$

(see Section 2.2 for a formal definition). Diaconis and Saloff-Coste [5] and Pak and Sidenko [12] use the dual characterization of $W_{d,1}$ [15], Remark 6.5, that is especially relevant for simulations:

$$(1) \quad W_{d,1}(\mu,\nu) = \sup\left\{\int_M f\,d(\mu-\nu) : f : M \to \mathbb{R} \text{ is } 1-\text{Lipschitz under } d\right\}.$$

That is, if one can sample from $\mu$, we can estimate $\int_M f\,d\nu$ for any Lipschitz $f$ up to a $W_{d,1}(\mu,\nu)$ intrinsic bias. This is a natural metric for many applications; as a case in point, we briefly discuss a suggestion of Ailon and Chazelle [1]. It is well known that one can "reduce the dimension" of a point set $S \subset \mathbb{R}^n$ while approximately preserving distances by first applying a random linear transformation $X$ drawn from the Haar measure on $SO(n)$ and then projecting onto the first $k$ coordinates. A result known as the Johnson Lindenstrauss lemma says that if one chooses $k = O(\ln|S|/\varepsilon^2)$ (which does not depend on the ambient dimension $n$), then the ratios of pairwise distances in $S$ are all preserved up to $(1\pm\varepsilon)$-factors, with high probability. One can easily check that a similar result holds when $X$ is sufficiently close to being Haar distributed in the $W_{d,1}$ metric (for an appropriate metric $d$; see below). As noted in [1], for $X = X(t)$ coming from Kac's random walk, the products $s_t = X(t)s$ ($s \in S$) can be computed with just a constant amount of extra memory, as the map $s_t \mapsto s_{t+1}$ affects only two coordinates of $s_t$; hence, if we can prove that $X(t)$ converges rapidly to a Haar measure in the $W_{d,1}$ distance, we have a time- and memory-efficient way of doing dimensionality reduction.

Our main result is a rapid mixing bound for the $SO(n)$ walk. We consider $M = SO(n)$ with two different choices of metric $d$. For $a,b \in SO(n)$ we define:

$$\text{hs}(a,b) \equiv \|a-b\|_{\text{hs}} = \sqrt{\text{Tr}((a-b)^\dagger(a-b))} \quad \text{the Hilbert–Schmidt norm;}$$

$D(a,b) \equiv$ the Riemannian metric on $SO(n)$ induced by the Hilbert–Schmidt
  inner product $\langle u,v\rangle_{\text{hs}} \equiv \text{Tr}(u^\dagger v)$.

Clearly, $\text{hs} \leq D$ always. Define the $L^p$ *transportation-cost mixing times*:

$$\tau_{d,p}(\varepsilon) \equiv \inf\{t \in \mathbb{N} : W_{d,p}(\mu K^t, \mathcal{H}) \leq \varepsilon \text{ for all prob. measures } \mu \text{ on } SO(n)\},$$

where $d = D$ or hs; $\mathcal{H}$ is the Haar measure on $SO(n)$; $K$ is the transition kernel for Kac's walk; and $\mu K^t$ is the time-$t$ distribution of a walk started from distribution $\mu$. Note that $\tau_{\text{hs},p}(\cdot) \leq \tau_{D,p}(\cdot)$ and that both mixing times are increasing in $p$. We will show the following:

THEOREM 1 (Main result). *For all $n \in \mathbb{N} \setminus \{0, 1\}$, Kac's random walk on $SO(n)$ satisfies the following mixing time estimate:*

$$\tau_{D,2}(\varepsilon) \leq \left\lceil n^2 \ln\left(\frac{\pi\sqrt{n}}{\varepsilon}\right) \right\rceil.$$

Thus, $O(n^2 \ln n)$ steps of the Markov chain suffice to bring $\mu K^t$ $\varepsilon$-close to the Haar measure $\mathcal{H}$ for any $\varepsilon = n^{-O(1)}$. This improves upon the $O(n^4 \ln n)$ bound by Diaconis and Saloff-Coste [5] and a very recent preprint by Pak and Sidenko [12] that lowered the estimate to $O(n^{2.5} \ln n)$ (we only learned about that result after proving the main results in the present paper). Moreover, these two papers treated only the $L^1$ case for $d = \text{hs}$, whereas we consider the stronger $L^2$ case with the stronger metric $D$.

We also show that our bounds are tight up to a $O(\ln n)$ factor, for all $n^{-O(1)} \leq \varepsilon \leq \varepsilon_0$ ($\varepsilon_0$ some constant), even when applied to $p = 1$ and $d = \text{hs}$.

THEOREM 2. *There exist $c, \varepsilon_0 > 0$ such that, for all $n \in \mathbb{N} \setminus \{0, 1\}$,*

$$\tau_{\text{hs},1}(\varepsilon_0) \geq cn^2.$$

Theorem 2 follows from a general lower bound for the mixing time of random walks induced by group actions. The general result might be already known, but since we could not find a proof of it elsewhere, we provide our own proof in Section 6. The bound in Theorem 2 was also claimed in [12].

The key to proving our main result, Theorem 1, is a contraction property of the Markov transition kernel of the random walk under consideration. Fix again a metric space $(M, d)$. For $\xi > 0$, say that a Markov transition kernel $P$ on $M$ is $\xi$-*Lipschitz for the $W_{d,p}$ metric* if for all probability measures $\mu, \nu$ on $M$ with finite $p$th moments (cf. Section 2.2)

(2) $$W_{d,p}(\mu P, \nu P) \leq \xi W_{d,p}(\mu, \nu).$$

If $\xi < 1$, we shall also say that $P$ is $\xi$-*contracting*. We will prove the following estimate:

LEMMA 1. *In the same setting as Theorem 1, Kac's random walk on matrices is*

$$\sqrt{1 - \frac{1}{\binom{n}{2}}}\text{-contracting}$$

*in the $W_{D,2}$ metric.*

The proof of Lemma 1 follows a strategy related to the *path coupling method* for discrete Markov chains introduced by Bubley and Dyer [3]. Suppose $P$ is now a Markov chain on the set of vertices $V$ of a connected graph $G$. The graph induces a natural shortest-path metric $d$ on $G$. It is sometimes possible to prove a "local contraction" estimate of the following form: for any $x, y \in V$ that are adjacent in $G$, there is a coupling of $X$ (distributed according to one step of $P$ from $x$) and $Y$ (distributed according to one step of $P$ from $y$) such that

$$\mathbb{E}[d(X,Y)] \leq \xi = \xi d(x,y) < 1.$$

If that is the case, Bubley and Dyer proved that the local couplings extend to "globally contracting" couplings for all random pairs $(x,y) = (X_0, Y_0) \in V^2$, with

$$\mathbb{E}[d(X,Y)] \leq \xi \mathbb{E}[d(X_0, Y_0)].$$

This implies, in particular, that $W_{d,1}(\mu P^t, \nu P^t) \leq \text{diam}(G)\xi^t$ for all distributions $\mu$, $\nu$, where $\text{diam}(G)$ is the diameter of the graph $G$. In the discrete setting such results easily extend to total variation bounds.

Our adaptation of their technique is based on the fact that $SO(n)$ is a *geodesic space* with the metric $D$: that is, $D(a,b)$ is the length of the shortest curve connecting $a$ and $b$. We will show that whenever $(M, d)$ is a geodesic space (or more generally a *length space*; see Section 2.1) and $P$ is such that, for all deterministic $x, y \in M$ with $d(x,y) \ll 1$,

$$\mathbb{E}[d(X,Y)^p] \leq (\xi + o(1))^p d(x,y)^p,$$

then $P$ is $\xi$-contracting and $W_{d,p}(\mu P^t, \eta P^t) \leq \xi^t \text{diam}(M)$ for all probability measures $\mu, \eta$ with finite $p$th moments, where $\text{diam}(M)$ is the diameter of $M$. That is, we show that if $(M,d)$ is a Polish length space and $P$ satisfies some reasonable assumptions, one can always extend "local contracting couplings" of $P$-steps from nearby deterministic states to "global contracting couplings" for arbitrary initial distributions. This result is stated as Theorem 3 below.

As with the original path-coupling methodology, proving local contraction is the problem-specific part of our technique. For Kac's walk, one can use the *local geometry* of $SO(n)$ as a Riemannian manifold to do calculations in the *tangent space*, which greatly simplifies our proof. The same idea can be applied to two related random walks (discussed in Section 5):

- a variant of Kac's walk where $\theta_t$ is nonuniform;
- a random walk on the set $U(n)$ of $n \times n$ unitary matrices where each step consists of applying a unitary transformation from $U(2)$ to the span of a pair of coordinate vectors.

Pak and Sidenko [12] use a related coupling construction, but neither use the local structure of $SO(n)$ as effectively, nor do they state any general result on local-to-global couplings. Diaconis and Saloff-Coste [5] use the analytic technique known as the *comparison method*, which seems intrinsically suboptimal for this problem, as well as more difficult to apply. [These two papers also handle some variants of Kac's process which do not seem to be related to the case we consider in Section 5.]

The general idea of contracting Markov chains with continuous state spaces has appeared in other works. Particularly relevant is a preprint of Ollivier [11], released during the preparation of the present paper, that contains a result related to (but a bit weaker than) our "path coupling" result, Theorem 3. That paper is devoted to the study of "positive Ricci curvature" for Markov chains on metric spaces, which is precisely what we call $\xi$-contractivity for $W_{d,1}$; from that one can deduce many properties, such as concentration for the stationary distribution and some log-Sobolev-like inequalities. See [11] for details and other references where contraction properties of Markov chains have been used recently. There have been many other recent results involving analytic, geometric and probabilistic applications of transportation cost [9, 13, 14]; this suggests that our techniques may find applications in that growing field. Of course, we also hope that our techniques will be applied to obtain mixing bounds of other Markov chains of intrinsic interest, not necessarily related to such geometric and analytic phenomena.

The remainder of the paper is as follows. Section 2 reviews some important concepts from probability, metric geometry and optimal transport. Section 3 proves our general result on local-to-global couplings, Theorem 3. Section 4 contains the definition of Kac's random walk on matrices and the proofs of Lemma 1 and Theorem 1. Section 5 sketches the two other random walks described above. Mixing time lower bounds are discussed in Section 6. Finally, Section 7 discusses other applications of our method and presents an open problem.

## 2. Preliminaries.

2.1. *Metric spaces, length spaces, $\sigma$-fields.* Whenever we discuss metric spaces $(M, d)$, saying that $A \subset M$ is *measurable* will mean that $A$ belongs to the $\sigma$-field generated by open sets in $M$, that is, the Borel $\sigma$-field $\mathcal{B}(M)$. Moreover, all measures on metric spaces will be implicitly defined over Borel sets. We will always assume that the metric spaces under consideration are *Polish*, that is, complete and separable.

Let $\gamma : [a, b] \to M$ be a continuous curve. The *length* $L_d(\gamma)$ of $\gamma$ (according to the metric $d$) is the following supremum:

$$L_d(\gamma) \equiv \sup\left\{\sum_{i=1}^n d(\gamma(t_{i-1}), \gamma(t_i)) : n \in \mathbb{N}, a = t_0 \leq t_1 \leq t_2 \leq \cdots \leq t_n = b\right\}.$$

The curve $\gamma$ is *rectifiable* if $L_d(\gamma) < +\infty$. The metric space $(M,d)$ is a *length space* if for all $x, y \in M$

$$d(x,y) = \inf\{L_d(\gamma) : \gamma : [0,1] \to M \text{ continuous}, \gamma(0) = x, \gamma(1) = y\}.$$

All complete Riemannian manifolds and their Gromov–Hausdorff limits are length spaces. Nonlocally-compact examples of Polish length spaces include separable Hilbert spaces, as well as infinite-dimensional $L_1$ spaces.

### 2.2. Distributions, couplings and mass transportation.
All facts stated below can be found in [15], Chapter 6.

Let $(M,d)$ be a metric space and $\Pr(M)$ be the space of probability measures on (the Borel $\sigma$-field of) $M$. Given $\mu, \nu \in \Pr(M)$, a measure $\nu \in \Pr(M \times M)$ (with the product Borel $\sigma$-field) is a *coupling* of $(\mu, \nu)$ if for all Borel-measurable $A \subset M$,

$$\eta(A \times M) = \mu(A), \qquad \eta(M \times A) = \nu(A).$$

The set of couplings of $(\mu, \nu)$ is denoted by $\mathrm{Cp}(\mu, \nu)$. This is always a nonempty set since the product measure $\mu \times \nu$ is in it.

Given $p \geq 1$, $\Pr_{d,p}(M) \subset \Pr(M)$ is the set of all probability measures $\mu$ with finite $p$th moments, that is, such that for some (and hence all) $o \in M$,

$$\int_M d(o,x)^p \, d\mu(x) < +\infty.$$

One can define the $L^p$ *transportation cost* (or $L^p$ *Wasserstein*) *metric* $W_{d,p}$ on $\Pr_{d,p}(M)$ by the formula

$$W_{d,p}(\mu,\nu)^p \equiv \inf\left\{ \int_{M \times M} d(x,y)^p \, d\eta(x,y) : \eta \in \mathrm{Cp}(\mu,\nu) \right\},$$
(3)
$$\mu, \nu \in \Pr_{d,p}(M).$$

Such metrics are related to the "mass transportation problem" where one attempts to minimize the average distance traveled by grains of sand when a sandpile is moved from one configuration to another.

It is known that $(\Pr_{d,p}(M), W_{d,p})$ is Polish iff $(M,d)$ is Polish. If $(M,d)$ is Polish, the infimum above is always achieved by some $\eta = \eta^{\mathrm{opt}}(\mu, \nu)$, which we will refer to as a $L^p$-*optimal coupling* of $\mu$ and $\nu$.

For $x \in M$, $\delta_x \in \Pr(M)$ is the *Dirac delta* (or *point mass*) at $x$, the distribution that assigns measure 1 to the set $\{x\}$. A basic property of mass transportation is that if $x, y \in M$, then

$$W_{d,p}(\delta_x, \delta_y) = d(x,y).$$

If $\mu \in \Pr_{d,p}(M)$ and $\delta_x$ is as above,

$$W_{d,p}(\delta_x, \mu)^p = \int d(x,y)^p \, d\mu(y).$$

It is often convenient to deal with random variables rather than measures. If $X$ is a $M$-valued random variable,

$$\mathcal{L}_X \in \Pr(M)$$

is the distribution (or law) of $X$. Notice that

$$\mathcal{L}_X \in \Pr_{d,p}(M) \quad \Leftrightarrow \quad \mathbb{E}[d(o,X)^p] < +\infty \quad \text{for some (all) } o \in M.$$

We will write

$$X =_{\mathcal{L}} \mu$$

whenever $X$ is a random variable with $\mathcal{L}_X = \mu$. Call a random pair $(X,Y)$ a *coupling* of $(\mu,\nu)$ if $\mathcal{L}_{(X,Y)} \in \mathrm{Cp}(\mu,\nu)$. $W_{d,p}(\mu,\nu)$ can be equivalently viewed as the infimum of $\mathbb{E}[d(X,Y)^p]^{1/p}$ over all such couplings.

Finally, we note that if $M$ is compact (as it is in our main application), then for any $p \geq 1$ $\Pr_{d,p}(M) = \Pr(M)$ and $W_{d,p}$ metrizes weak convergence.

2.3. *Markov transition kernels.* In this section we assume $(M,d)$ is Polish. A *Markov transition kernel* on $M$ is a map

$$P : M \times \mathcal{B}(M) \to [0,1]$$

such that, for all $x \in M$, $P_x(\cdot) \equiv P(x,\cdot)$ is a probability measure and for all $A \in \mathcal{B}(M)$, $P_x(A)$ is a measurable function of $x$. A Markov transition kernel defines an $M$-valued Markov chain: for each $\mu \in \Pr(M)$, there exists a unique distribution on sequences of $M$-valued random variables

$$\{X(t)\}_{t=0}^{+\infty}$$

such that $X(0) =_{\mathcal{L}} \mu$ and for all $t \in \{1,2,3,\ldots\}$, the distribution of $X(t)$ conditioned on $\{X(s)\}_{s=0}^{t-1}$ is $P_{X(t-1)}$.

For $\mu \in \Pr(M)$ and $t \in \mathbb{N}$, $\mu P^t$ is the measure of $X(t)$ defined as above; one can check that $\mu P^{t+1} = (\mu P^t)P$ for all $t \geq 0$.

**3. From local to global couplings.** In this section we will discuss our method for moving from local to global bounds for the Lipschitz properties of Markov kernels. In our application we have a Markov kernel $P$ on a Polish space $(M,d)$. Using explicit couplings, we will show that, for some $C > 0$ and all $x,y \in M$,

$$W_{d,p}(P_x, P_y) \leq (C + o(1))d(x,y),$$

where $o(1) \to 0$ when $y \to x$. The main result in this section implies that, under some natural conditions, it follows that $W_{d,p}(\mu P, \nu P) \leq Cr$ whenever $\mu, \nu \in \Pr_{d,p}(M)$ are $r$-close.

We first state a definition.

DEFINITION 1. A map $f:M \to N$ between metric spaces $(M,d)$ and $(N,d')$ is said to be *locally C-Lipschitz* (for some $C > 0$) if for all $x \in M$

$$\limsup_{y \to x} \frac{d'(f(x), f(y))}{d(x,y)} \leq C.$$

THEOREM 3 (Local-to-global coupling). *Suppose $(M,d)$ is a Polish length space, $p \geq 1$ is given and $P$ is a Markov transition kernel on $(M,d)$ satisfying the following characteristics:*

1. *$P_x$ has finite pth moments for all $x$: that is, $P_x \in \Pr_{d,p}(M)$ for all $x \in M$;*
2. *$P$ is locally $C$-Lipschitz on $M$. That is, the map*

$$x \mapsto P_x$$

*from $(M,d)$ to $(\Pr_{d,p}(M), W_{d,p})$ is locally $C$-Lipschitz.*

*Then for all $\mu \in \Pr_{d,p}(M)$, we also have $\mu P \in \Pr_{d,p}(M)$ and, moreover, the map $\mu \mapsto \mu P$ is $C$-Lipschitz, that is,*

$$\forall \mu, \nu \in \Pr_{d,p}(M), \qquad W_{d,p}(\mu P, \nu P) \leq C W_{d,p}(\mu, \nu).$$

Before we prove this result, we discuss its application to the setting where $C = (1-\kappa)$ for some $\kappa > 0$, the diameter $\text{diam}_d(M)$ of $(M,d)$ is bounded (Ollivier [11] noted that, for $C = (1-\kappa) < 1$, $\text{diam}_d(M) \leq 2\Delta/\kappa$, where $\Delta = \sup_{x \in M} W_{d,p}(\delta_x, P_x)$. Hence, the assumption that $\text{diam}_d(M) < +\infty$ is equivalent to $\Delta < +\infty$) and the other assumptions of Theorem 3 are satisfied. In this case $\Pr_{d,p}(M) = \Pr(M)$, that is, all probability measures have bounded $p$th moments. Moreover, Banach's fixed point theorem states that a $(1-\kappa)$-Lipschitz map from a complete metric space to itself has a unique fixed point. Since $(\Pr(M), W_{d,p})$ is Polish and $\mu \mapsto \mu P$ is a $(1-\kappa)$-Lipschitz map from this space to itself, there exists a unique element $\mu_* \in \Pr(M)$ with $\mu_* P = \mu_*$.

It follows that $\mu_*$ is the unique $P$-invariant distribution on $M$. Moreover, for all $t \in \mathbb{N}$ and $\mu \in \Pr(M)$,

$$W_{d,p}(\mu P^t, \mu_*) = W_{d,p}(\mu P^t, \mu_* P^t) \leq (1-\kappa)^t W_{d,p}(\mu, \mu_*) \leq (\text{diam}_d(M)) e^{-\kappa t}.$$

We collect those facts in the following corollary.

COROLLARY 1. *Assume $(M,d)$ and $P$ satisfy the assumptions of Theorem 3 for some $p \geq 1$ and $C = (1-\kappa) < 1$ (i.e., $\kappa > 0$). Assume, moreover, that the diameter $\text{diam}_d(M)$ of $(M,d)$ is finite. Then there exists a unique $P$-invariant measure $\mu_*$ on $M$. Moreover, the $L^p$ transportation cost mixing times*

$$\tau_{d,p}(\varepsilon) \equiv \min\{t \in \mathbb{N} : \forall \mu \in \Pr(M), W_{d,p}(\mu P^t, \mu_*) \leq \varepsilon\}$$

*satisfy*

$$\tau_{d,p}(\varepsilon) \leq \left\lceil \kappa^{-1} \ln\left(\frac{\text{diam}_d(M)}{\varepsilon}\right) \right\rceil.$$

We now proceed to prove the theorem.

PROOF OF THEOREM 3. The first step of the proof is a simple lemma (proven subsequently) about locally Lipschitz functions.

LEMMA 2. *With the notation of Definition 1, assume that $M$ is a length space. Then any $f : M \to N$ that is locally $C$-Lipschitz is $C$-Lipschitz according to the standard definition.*

For our proof we only need the following direct consequence [let $(N, d') = (\text{Pr}_{d,p}(M), W_{d,p})$, $f(x) = P_x$].

COROLLARY 2. *If $P$ is a Markov transition kernel on a length space $(M, d)$ satisfying condition 2 of Theorem 3, then $W_{d,p}(P_x, P_y) \leq Cd(x, y)$ for all $x, y \in M$.*

The bounding of $W_{d,p}(P_x, P_y)$ can be thought of as an implicit construction of a coupling along a geodesic path; this is precisely where the name "path coupling" comes from.

The second lemma we need (proven in Section 3.2) shows that $\mu P \in \text{Pr}_{d,p}(M)$ whenever $\mu \in \text{Pr}_{d,p}(M)$ and shows that we will only need to compare $\mu P$ and $\nu P$, for $\mu$, $\nu$ with countable support.

LEMMA 3. *Let $(M, d)$ be Polish. Suppose $P$ is a Markov transition kernel on $M$ such that:*

1. *$P_x \in \text{Pr}_{d,p}(M)$ for all $x \in M$;*
2. *$x \mapsto P_x$ is a $C$-Lipschitz map from $M$ to $\text{Pr}_{d,p}(M)$.*

*Then for all $\mu \in \text{Pr}_{d,p}(M)$ we have $\mu P \in \text{Pr}_{d,p}(M)$. Moreover, there exists a sequence $\{\mu_j\}_j \subset \text{Pr}_{d,p}(M)$ of measures with countable support such that $W_{d,p}(\mu_j, \mu) \to 0$ and $W_{d,p}(\mu_j P, \mu P) \to 0$.*

The lemma implies the following statement: if $W_{d,p}(\mu P, \nu P) \leq CW_{d,p}(\mu, \nu)$ for all $\mu$, $\nu$ in $\text{Pr}_{d,p}(M)$ that have countable support, then the same holds for all $\mu$, $\nu$ in $\text{Pr}_{d,p}(M)$. Our final goal is to prove the Lipschitz estimate for measures with countable support.

Thus, let $\mu = \sum_{j \in \mathbb{N}} p_j \delta_{x_j}$ be a convex combination of a countable number of point masses ($x_j \in M$ for all $j$); similarly, let $\nu = \sum_{k \in \mathbb{N}} q_k \delta_{y_k}$. The $L^p$-optimal coupling $\eta$ of $\mu$ and $\nu$ is of the form

$$\eta = \sum_{j,k \in \mathbb{N}} r_{j,k} \delta_{(x_j, y_k)}$$

for some convex weights $r_{j,k}$. Now define for each pair $j, k$ a $L^p$-optimal coupling $Q_{j,k}$ of $P_{x_j}, P_{y_k}$. Then

$$\eta' = \sum_{j,k \in \mathbb{N}} r_{j,k} Q_{j,k} \in \mathrm{Cp}(\mu P, \nu P).$$

Moreover, since $x \mapsto P_x$ is $C$-Lipschitz,

$$\int_{M \times M} d(u,v)^p \, dQ_{j,k}(u,v) = W_{d,p}(P_{x_j}, P_{y_k})^p \leq C^p d(x_j, y_k)^p,$$

which implies

$$\begin{aligned}
W_{d,p}(\mu P, \nu P)^p &\leq \int_{M \times M} d(u,v)^p \, d\eta'(u,v) \\
&= \sum_{j,k \in \mathbb{N}} r_{j,k} \int_{M \times M} d(u,v)^p \, dQ_{j,k}(u,v) \\
&\leq C^p \sum_{j,k \in \mathbb{N}} r_{j,k} d(x_j, y_k)^p \\
&= C^p \int_{M \times M} d(u,v)^p \, d\eta(u,v).
\end{aligned}$$

The RHS is simply $C^p W_{d,p}(\mu, \nu)^p$. □

REMARK 1. Ollivier presents a similar result for $p = 1$ in [11], Proposition 17. His proof relies on a quite nontrivial fact (proven in, e.g., [15]): the existence of a Markov transition kernel $Q$ on $M^2$ such that, for all $(x, y) \in M^2$, $Q_{(x,y)}$ is a 1-optimal coupling of $(P_x, P_y)$. Our argument provides an alternative approach, which is perhaps simpler, to the same result. Moreover, his proposition implies our theorem only when $P$ satisfies:

$$\limsup_{r \searrow 0} \sup_{x,y \in M : d(x,y) \leq r} \frac{W_{d,1}(P_x, P_y)}{r} \leq C,$$

which is a stronger requirement than our local Lipschitz condition.

3.1. *Proof of Lemma 2.* It suffices to show that, for all $x, y \in M$, any continuous curve $\gamma \colon [0,1] \to M$ connecting $\gamma(0) = x$ to $\gamma(1) = y$ and any number $C' > C$,

$$W_{d,p}(P_x, P_y) \leq C' L_d(\gamma).$$

To prove this, assume without loss of generality that $L_d(\gamma) < +\infty$. For $0 \leq t_1 < t_2 \leq 1$, define the *length function*

$$\ell(t_1, t_2) \equiv L_d(\gamma \mid_{[\min\{t_1,t_2\}, \max\{t_1,t_2\}]}).$$

It is an exercise to show that $\ell$ is a continuous function $\ell(t_1, t_2) \geq d(\gamma(t_1), \gamma(t_2))$ and

(4) $\quad \forall 0 \leq t_1 < t_2 < t_3 \leq 1, \qquad \ell(t_1, t_2) + d(\gamma(t_2), \gamma(t_3)) \leq \ell(t_1, t_3).$

For each $t \in [0,1]$, we have

$$\limsup_{s \to t} \frac{d'(f(\gamma(s)), f(\gamma(t)))}{\ell(s,t)} \leq \limsup_{s \to t} \frac{d'(f(\gamma(s)), f(\gamma(t)))}{d(\gamma(s), \gamma(t))} \leq C,$$

by the local Lipschitz assumption. Since $C' > C$, one can find, for any $t \in [0,1)$, some $\delta_t \in (0, 1-t)$ such that $\forall s \in (t, t+\delta_t]$, $d'(f(\gamma(s)), f(\gamma(t)))) \leq C'\ell(t,s)$.

Now set

$$T \equiv \sup\{t \in [0,1] \colon d'(f(\gamma(0)), f(\gamma(t))) \leq C'\ell(0,t)\}.$$

Notice that

(5) $\qquad d'(f(\gamma(0)), f(\gamma(T))) \leq C'\ell(0,T)$

by continuity. We claim that $T = 1$. To see this, suppose $T < 1$ and set $\delta = \delta_T$. Then

$$d'(f(\gamma(0)), f(\gamma(T+\delta))) \leq d'(f(\gamma(0)), f(\gamma(T)))$$
$$+ d'(f(\gamma(T)), f(\gamma(T+\delta)))$$

[use (5) and defn. of $\delta_T$] $\leq C'\ell(0,T) + C'd(\gamma(T), \gamma(T+\delta))$

[use (4)] $= C'\ell(0, T+\delta),$

which contradicts the fact that $T$ is the supremum of the corresponding set. We deduce that $T = 1$ and

$$d'(f(x), f(y)) = d'(f(\gamma(0)), f(\gamma(1))) \leq C'\ell(0,1) = C'L_d(\gamma),$$

as desired.

3.2. *Proof of Lemma 3.* For the first statement we note that, for a given reference point $y \in M$,

$$A^p \equiv \int_M d(y,z)^p \, dP_y(z) < +\infty.$$

Now for any $x \in M$, let $(X,Y)$ be a $L^p$-optimal coupling of $(P_x, P_y)$. Then

$$\|d(y,X)\|_{L_p} \leq \|d(y,Y)\|_{L_p} + \|d(Y,X)\|_{L_p} = A + W_{d,p}(P_x, P_y) \leq A + Cd(x,y),$$

which is the same as

$$\int_M d(y,v)^p \, dP_x(v) \leq (A + Cd(x,y))^p.$$

Hence, if $\mu \in \mathrm{Pr}_{d,p}(M)$,

$$\int_M d(y,v)^p \, d\mu P(v) = \int_M \left( \int_M d(y,v)^p \, dP_x(v) \right) d\mu(x)$$

$$= \int_M \left( \int_M [A + Cd(y,v)]^p \, dP_x(v) \right) d\mu(x)$$

$$[\text{use } |a+b|^p \leq 2^p(|a|^p + |b|^p)] \leq (2C)^p \int_M \int_M d(y,v)^p \, dP_x(v) \, d\mu(x)$$

$$+ 2^p \int_M \int_M A^p \, dP_x(v) \, d\mu(x)$$

$$\leq (2C)^p \int_M d(y,v)^p \, d\mu(x) + (2A)^p$$

$$[\mu \in \mathrm{Pr}_{d,p}(M)] < +\infty.$$

Thus, $\mu P$ is in $\mathrm{Pr}_{d,p}(M)$ whenever $\mu$ is.

We now present a discrete approximation scheme for $\mu$ and $\mu P$. Since $M$ is separable, there exists a sequence of partitions $\{\mathcal{P}_j\}_{j \in \mathbb{N}}$ of $M$ such that:

- each partition contains countably many measurable sets;
- for all $j \in \mathbb{N}$, $\mathcal{P}_{j+1}$ refines $\mathcal{P}_j$; and
- for all $j \in \mathbb{N}$, the sets in $\mathcal{P}_j$ have diameter at most $\varepsilon_j$ for some sequence $\varepsilon_j \to 0$.

Let us also assume that for each $j \in \mathbb{N}$ and $A \in \mathcal{P}_j$ we have picked some $x_A^{(j)} \in A$. Consider the measures

(6) $$\mu_j \equiv \sum_{A \in \mathcal{P}_j} \mu(A) \delta_{x_A^{(j)}}.$$

Clearly, $\mu_j \in \mathrm{Pr}_{d,p}(M)$ for all $j$ and $W_{d,p}(\mu_j, \mu) \to 0$ when $j \to +\infty$. Our goal will be to show that $W_{d,p}(\mu_j P, \mu P) \to 0$. First recall that $x \mapsto P_x$ is $C$-Lipschitz, hence, if $x, y \in M$ and $d(x,y) \leq \varepsilon_j$, $W_{d,p}(P_x, P_y) \leq C\varepsilon_j$. In particular, for all $j \in \mathbb{N}$, all $A \in \mathcal{P}_j$ and all $x \in A$,

$$W_{d,p}(P_{x_A^{(j)}}, P_x) \leq C\varepsilon_j.$$

We will use this to show that

(7)
$$\forall j < k, \qquad W_{d,p}(\mu_j P, \mu_k P) \leq C\varepsilon_j$$

(in particular, $\{\mu_j P\}_j$ is Cauchy).

Recall that if $j < k$, $\mathcal{P}_k$ is a refinement of $\mathcal{P}_j$, hence, for all $B \in \mathcal{P}_k$ there exists a set $A_B \in \mathcal{P}_j$ with $B \subset A_B$. For each such $B$, we have $x_B^{(k)} \in A_B$, which has diameter $\leq \varepsilon_j$, hence, $d(x_B^{(k)}, x_{A_B}^{(j)}) \leq \varepsilon_j$ and there exists a coupling $\eta_{B,k,j}$ of $P_{x_B^{(k)}}$ and $P_{x_{A_B}^{(j)}}$ with

$$\int_{M \times M} d(u,v)^p \, d\eta_{B,k,j}(u,v) = W_{d,p}(P_{x_B^{(k)}}, P_{x_{A_B}^{(j)}})^p \leq (C\varepsilon_j)^p.$$

Extend this to a coupling of $\mu_k P$ and $\mu_j P$ by

$$\eta_{k,j} \equiv \sum_{B \in \mathcal{P}_k} \mu(B) \eta_{B,k,j}.$$

To prove that $\eta_{k,j} \in \mathrm{Cp}(\mu_j P, \mu_k P)$, notice that the first marginal of this measure is

$$\sum_{B \in \mathcal{P}_k} \mu(B) P_{x_B^{(k)}} = \mu_k P.$$

Moreover, for any $A \in \mathcal{P}_j$, the set of all $B \in \mathcal{P}_k$ with $A_B = A$ is a partition of $A$, hence the second marginal is also right:

$$\sum_{B \in \mathcal{P}_k} \mu(B) P_{x_{A_B}^{(j)}} = \sum_{A \in \mathcal{P}_j} \left( \sum_{B \in \mathcal{P}_k : A_B = A} \mu(B) \right) P_{x_A^{(j)}} = \sum_{A \in \mathcal{P}_j} \mu(A) P_{x_A^{(j)}} = \mu_j P.$$

It follows that $\eta_{k,j} \in \mathrm{Cp}(\mu_j P, \mu_k P)$ and, moreover, one can check that

$$\int_{M \times M} d(u,v)^p \, d\eta_{k,j}(u,v) \leq (C\varepsilon_j)^p,$$

which implies (7).

$(\mathrm{Pr}_{d,p}(M), W_{d,p})$ is Polish since $(M, d)$ is. By the above, we know that there exists a measure $\alpha \in \mathrm{Pr}_{d,p}(M)$ such that $W_{d,p}(\mu_j P, \alpha) \to 0$. This also implies [15], Theorem 6.8, that $\mu_j P \Rightarrow \alpha$ in the weak topology. However, it is an exercise to show that $\mu_j P \Rightarrow \mu P$ weakly, hence, $\alpha = \mu P$ and $W_{d,p}(\mu_j P, \mu P) \to 0$, as desired.

**4. Analysis of Kac's random walk.**

4.1. *Definitions.* Let $M(n,\mathbb{R})$ be the set of all $n \times n$ matrices with real-valued entries. These are the linear operators from $\mathbb{R}^n$ to itself and we equip $\mathbb{R}^n$ with a canonical basis $e_1, \ldots, e_n$ of orthonormal vectors. For $a \in M(n,\mathbb{R})$, $a^\dagger$ is the transpose of $a$ in the basis $e_1, \ldots, e_n$. Using it, one can define the Hilbert–Schmidt inner product $\langle a, b \rangle_{\mathrm{hs}} \equiv \mathrm{Tr}(a^\dagger b)$ on $M(n,\mathbb{R})$, under which it is isomorphic to $\mathbb{R}^{n^2}$ with the standard Euclidean inner product. We let $\|\cdot\|_{\mathrm{hs}}$ be the corresponding norm.

An element $a \in M(n,\mathbb{R})$ is *orthogonal* if $aa^\dagger = \mathrm{id}$, the identity matrix. The subset of $M(n,\mathbb{R})$ given by

$$SO(n) \equiv \{a \in M(n,\mathbb{R}) : aa^\dagger = \mathrm{id}, \det(a) = 1\}$$

is a smooth, compact, connected submanifold of $M(n,\mathbb{R})$. It is also a Lie group since it is closed under matrix multiplication and matrix inverse. Therefore, $SO(n)$ has a Haar measure $\mathcal{H}$, which we may define as the unique probability measure on that group such that, for all measurable $S \subset SO(n)$ and $a \in SO(n)$, we have $\mathcal{H}(S) = \mathcal{H}(Sa) = \mathcal{H}(aS)$, where $Sa = \{sa : s \in S\}$ and $aS = \{as : s \in S\}$.

We now define *Kac's random walk on $SO(n)$*. For $1 \leq i < j \leq n$ and $\theta \in [0, 2\pi]$ define $R(i,j,\theta)$ as a rotation by $\theta$ of the plane generated by $e_i, e_j$. This is equivalent to setting

$$(8) \qquad R(i,j,\theta)e_k = \begin{cases} \cos\theta\, e_i + \sin\theta\, e_j, & k = i, \\ \cos\theta\, e_j - \sin\theta\, e_i, & k = j, \\ e_k, & k \in \{1, \ldots, n\} \setminus \{i, j\}, \end{cases}$$

and extending $R(i,j,\theta)$ to all $\psi \in \mathbb{R}^n$ by linearity. Kac's random walk on matrices corresponds to the following Markov transition kernel:

$$K_x(S) \equiv \frac{1}{2\pi \binom{n}{2}} \sum_{1 \leq i < j \leq n} \int_0^{2\pi} \delta_{R(i,j,\theta)x}(S)\, d\theta$$

$$[x \in SO(n), S \subset SO(n) \text{ measurable}].$$

Thus, to generate $X =_{\mathcal{L}} K_x$ from $x$, one chooses $1 \leq i < j \leq n$ uniformly at random from all $\binom{n}{2}$ possible choices, then picks $\theta \in [0, 2\pi]$ also uniformly at random and then sets $X = R(i,j,\theta)x$. The required measurability conditions are easily established. One can also check that the Haar measure $\mathcal{H}$ is $K$-invariant.

4.2. *The geometry of $SO(n)$.* We collect some standard facts that will be used in our proofs.

The tangent space at the identity matrix id is the set of all anti-self-adjoint operators

$$(9) \qquad T \equiv T_{\mathrm{id}} SO(n) = \{h \in M(n, \mathbb{R}) : h^\dagger = -h\}.$$

We let $D$ be the Riemannian metric on $SO(n)$ induced by $\langle \cdot, \cdot \rangle_{\mathrm{hs}}$. Since $SO(n)$ is compact, one can show the following:

$$(10) \quad \forall z, w \in SO(n), \qquad \|z - w\|_{\mathrm{hs}} \leq D(z, w) \leq \|z - w\|_{\mathrm{hs}} + O(\|z - w\|_{\mathrm{hs}}^2),$$

where $O(r^\alpha)$ is just some term whose absolute value is uniformly bounded by $c|r|^\alpha$ and $c > 0$ a constant not depending on $|r|$ (we will use this notation from now on). Moreover, if we let $\Pi_T$ be the orthogonal projector onto $T$ (according to the Hilbert–Schmidt inner product), then (although we will not use this fact, one can check that $\Pi_T \, \mathrm{id} = 0$)

$$(11) \qquad \forall z \in SO(n), \qquad \|z - \mathrm{id} - \Pi_T(z - \mathrm{id})\|_{\mathrm{hs}} \leq O(D(z, \mathrm{id})^2).$$

This is so because if $\|z - \mathrm{id}\|_{\mathrm{hs}} = r \ll 1$, then $\|z - \mathrm{id} - \tilde{h}\|_{\mathrm{hs}} = O(r^2)$ for some $\tilde{h} \in T$, and $\tilde{h} = h = \Pi_T(z - \mathrm{id})$ is the best choice of approximation one may make. Notice that the two equations together imply

$$(12) \qquad |D(z, \mathrm{id}) - \|\Pi_T(z - \mathrm{id})\|_{\mathrm{hs}}| = O(\|z - \mathrm{id}\|_{\mathrm{hs}}^2).$$

We notice that these distances are all invariant under multiplication: if $a, b, c \in SO(n)$,

$$D(ca, cb) = D(ac, bc) = D(a, b)$$

and similarly for $\mathrm{hs}(a, b) = \|a - b\|_{\mathrm{hs}}$.

### 4.3. *The contraction coefficient.* In this section we prove Lemma 1.

PROOF OF LEMMA 1. Consider $x, y \in SO(n)$ and let $D(x, y) = r$. Our main task is to show that there exists a coupling $(X, Y)$ of $(K_x, K_y)$ with

$$\mathbb{E}[D(X, Y)^2] \leq \left(1 - \frac{1}{\binom{n}{2}}\right) r^2 + O(r^3),$$

where, as in the previous section, $O(r^3)$ is some term that is uniformly bounded by a multiple of $|r|^3$. The existence of such a coupling implies that

$$W_{D,2}(K_x, K_y) \leq \sqrt{1 - \frac{1}{\binom{n}{2}}} D(x, y) + O(D(x, y)^2),$$

which shows that $K$ is locally $\sqrt{1 - 1/\binom{n}{2}}$-Lipschitz for $p = 2$.

Our coupling will be as follows. Suppose we set $X = R(i, j, \theta)x$ with $i, j, \theta$ randomly picked as prescribed in the definition of the random walk. We will set $Y = R(i, j, \theta')y$ with the same $i, j$ and some $\theta' = (\theta - \alpha) \bmod 2\pi$, where $\alpha = \alpha(i, j, x, y)$ depends on $i, j, x, y$ but *not* on $\theta$. In that case $\theta'$ is uniform on $[0, 2\pi]$ independently of $i, j, x, y$, hence, $(X, Y)$ is a valid coupling of $(K_x, K_y)$. Also notice that, using the invariance of $D$ under multiplication,

$$\begin{aligned}
D(X, Y) &= D(R(i, j, \theta)x, R(i, j, \theta')y) \\
&= D(R(i, j, \theta), R(i, j, \theta')yx^\dagger) = D(R(i, j, \alpha), yx^\dagger),
\end{aligned} \tag{13}$$

as

$$R(i, j, \theta')^\dagger R(i, j, \theta) = R(i, j, \theta - \theta') = R(i, j, \alpha).$$

We will use (10), (11) and (12) to bound the RHS of (13): this will allow us to do all calculations we need in the tangent space $T = T_{\mathrm{id}}SO(n)$. First, however, we need an orthonormal basis for that space. For each $1 \leq k < \ell \leq n$, let $a_{k\ell} \in T$ be the linear operator that is uniquely defined by

$$a_{k\ell}e_t \equiv \begin{cases} \dfrac{e_\ell}{\sqrt{2}}, & t = k, \\ -\dfrac{e_k}{\sqrt{2}}, & t = \ell, \\ 0, & t \in \{1, \ldots, n\} \setminus \{k, \ell\}. \end{cases}$$

One can check that $\{a_{k\ell}\}_{1 \leq k < \ell \leq n}$ is indeed an orthonormal basis for $T = T_{\mathrm{id}}SO(n)$ with the Hilbert–Schmidt inner product. For $1 \leq t \leq n$ we also define $d_t \in M(n, \mathbb{R})$ as the matrix that has a 1 at the $(t, t)$th entry and zeroes elsewhere. Then $\langle d_t, d_s \rangle_{\mathrm{hs}} = 1$ if $t = s$ and 0 otherwise and also $\langle d_t, a_{k\ell} \rangle_{\mathrm{hs}} = 0$ for any $t, k, \ell$. With these definitions, one can write

$$R(i, j, \alpha) = \mathrm{id} + (\cos\alpha - 1)d_i + (\cos\alpha - 1)d_j + \sqrt{2}\sin\alpha\, a_{ij}. \tag{14}$$

Now set $h = \Pi_T(yx^\dagger - \mathrm{id})$. Since $D(yx^\dagger, \mathrm{id}) = D(x, y) = r$, $\|h\|_{\mathrm{hs}} = r + O(r^2)$ and $\|yx^\dagger - \mathrm{id} - h\|_{\mathrm{hs}} = O(r^2)$. Suppose we commit ourselves to making a choice of $\alpha = O(r)$ (i.e., $|\alpha| \leq cr$ for a constant $c$ independent of $r$). Expanding sin and cos, we get

$$\|R(i, j, \alpha) - \mathrm{id} - \sqrt{2}\alpha a_{ij}\|_{\mathrm{hs}} = O(r^2).$$

Moreover, we also have

$$\begin{aligned}
D(yx^\dagger, R(i, j, \alpha)) & \\
&= \|yx^\dagger - R(i, j, \alpha)\|_{\mathrm{hs}} + O(\|yx^\dagger - R(i, j, \alpha)\|_{\mathrm{hs}}^2) \tag{15}\\
&= \|yx^\dagger - \mathrm{id} - \sqrt{2}\alpha a_{ij}\|_{\mathrm{hs}} + O(\|yx^\dagger - \mathrm{id} - \sqrt{2}\alpha a_{ij}\|_{\mathrm{hs}}^2 + r^2) \tag{16}\\
&= \|h - \sqrt{2}\alpha a_{ij}\|_{\mathrm{hs}} + O(r^2). \tag{17}
\end{aligned}$$

Thus, we choose $\alpha = \langle h, a_{ij} \rangle_{\rm hs}/\sqrt{2}$, which minimizes $\|h - \sqrt{2}\alpha a_{ij}\|_{\rm hs}$ and only depends on $i,j$ and $h = \Pi_T(yx^\dagger - {\rm id})$. Since the $a_{k\ell}$ form an orthonormal basis of $T \ni h$, we have

$$h = \sum_{1 \leq k < \ell \leq n} \langle h, a_{k\ell} \rangle_{\rm hs} a_{k\ell} \quad \Rightarrow \quad \sum_{1 \leq k < \ell \leq n} \langle h, a_{k\ell} \rangle_{\rm hs}^2 = \|h\|_{\rm hs}^2 = r^2 + O(r^3).$$

This shows that $|\alpha| = O(r)$ as desired and, moreover,

$$D(X,Y)^2 = D(yx^\dagger, R(i,j,\alpha))^2 \qquad [\text{by (13)}]$$
$$= \|h - \langle h, a_{ij} \rangle_{\rm hs} a_{ij}\|_{\rm hs}^2 + O(r^3)$$
$$(\text{expand } h) = \|h\|_{\rm hs}^2 - \langle h, a_{ij} \rangle_{\rm hs}^2 + O(r^3).$$

If we now average over $i,j,\theta$, we obtain

$$\mathbb{E}[D(X,Y)^2] = \|h\|_{\rm hs}^2 - \frac{1}{\binom{n}{2}} \sum_{1 \leq i < j \leq n} \langle h, a_{ij} \rangle_{\rm hs}^2 + O(r^3)$$
$$= \left(1 - \frac{1}{\binom{n}{2}}\right) \|h\|_{\rm hs}^2 + O(r^3)$$
$$= \left(1 - \frac{1}{\binom{n}{2}}\right) r^2 + O(r^3),$$

which is the desired bound.

To finish the proof, we apply our result on local-to-global couplings, Theorem 3. We have shown that the Markov transition kernel $P = K$ for Kac's random walk is locally $C$-Lipschitz for

$$C = \left(1 - \frac{1}{\binom{n}{2}}\right)^{1/2}, \qquad 1 \leq p \leq 2.$$

The remaining assumptions of Theorem 3 are trivially verified since $SO(n)$ has a bounded diameter. We conclude that

$$\forall \mu, \eta \in \Pr(SO(n)), \qquad W_{D,p}(\mu K, \nu K) \leq \sqrt{1 - \frac{1}{\binom{n}{2}}} W_{D,p}(\mu, \nu). \qquad \square$$

4.4. *Mixing time upper bound.* We now prove Theorem 1.

PROOF OF THEOREM 1. We shall apply Corollary 1 with $M = SO(n)$, $d = D$ and $P = K$. According to Lemma 1, we can take $C = \sqrt{1 - 1/\binom{n}{2}} \leq (1 - \kappa)$ for $\kappa = 1/n^2$.

We need an estimate for the diameter of $SO(n)$ under $D$. Let $a, b \in SO(n)$. Then $D(a,b) = D(c, {\rm id})$ with $c = ab^\dagger \in SO(n)$. It is well known that any such $c$ is a product of two-dimensional rotations on orthogonal subspaces; that is equivalent to saying that (after a change of basis of $\mathbb{R}^n$) one can write

$$c = \prod_{i=1}^{k} R(2i-1, 2i, \theta_i)$$

for $k = \lfloor n/2 \rfloor$ and $-\pi \leq \theta_i \leq \pi$ without loss of generality. Notice that one can rewrite this as [cf. (14)]

$$c = \sum_{i=1}^{k} [\cos\theta_i (d_{2i-1} + d_{2i}) + \sin\theta_i a_{2i-1,2i}].$$

Thus, the curve

$$\gamma(t) \equiv \sum_{i=1}^{k} [\cos t\theta_i (d_{2i-1} + d_{2i}) + \sin t\theta_i a_{2i-1,2i}], \qquad 0 \leq t \leq 1,$$

connects id to $c$ in $SO(n)$. Moreover, for all $0 \leq t \leq 1$,

$$\gamma'(t) \equiv \sum_{i=1}^{k} \theta_i [\cos(t\theta_i + \pi/2)(d_{2i-1} + d_{2i}) + \sin(t\theta_i + \pi/2) a_{2i-1,2i}]$$

and one can easily see that

$$\|\gamma'(t)\|_{\mathrm{hs}}^2 = 2\sum_{i=1}^{k} |\theta_i|^2 = 2k\pi^2 \leq \pi^2 n \qquad (\text{since } k \leq n/2).$$

We deduce that

$$\forall a, b \in SO(n), \qquad D(a,b) = D(ab^\dagger, \mathrm{id}) \leq \int_0^1 \|\gamma'(t)\|_{\mathrm{hs}}\, dt \leq \pi\sqrt{n}.$$

Thus, $\mathrm{diam}_D(SO(n)) \leq \pi\sqrt{n}$ and we deduce from the corollary that

$$\tau_{D,2}(\varepsilon) \leq \left\lceil n^2 \ln\left(\frac{\pi\sqrt{n}}{\varepsilon}\right) \right\rceil. \qquad \square$$

**5. Mixing bounds for other random walks.** In this section we briefly discuss the two random walks related to Kac's random walk mentioned in the introduction. Both proofs follow the previous one very closely and will be only sketched.

5.1. *Kac's walk with nonuniform angles.* Recall the definitions in Section 4.1. In this section we let $\rho : [0, 2\pi] \to \mathbb{R}^+$ be a density and define a variant $K^{(\rho)}$ of Kac's random walk on $SO(n)$ as follows:

$$K_x^{(\rho)}(S) \equiv \frac{1}{\binom{n}{2}} \sum_{1 \leq i < j \leq n} \int_0^{2\pi} \delta_{R(i,j,\theta)x}(S) \rho(\theta)\, d\theta$$

$$[x \in SO(n), S \subset SO(n) \text{ measurable}].$$

$K^{(\rho)}$ corresponds to picking the rotation angle with density $\rho$. One can check that $K^{(\rho)}$ is a valid Markov transition kernel for any density $\rho$ and that the original process corresponds to $\rho \equiv 1/2\pi$. We will prove the following:

THEOREM 4. *Suppose*

$$\rho_{\min} \equiv \min_{\theta \in [0,2\pi]} \rho(\theta) > 0.$$

*Then the $L^2$ transportation cost mixing time of $K^{(\rho)}$ satisfies*

$$\tau_{D,2}(\varepsilon) \leq \left\lceil \frac{n^2}{2\pi\rho_{\min}} \ln\left(\frac{\pi\sqrt{n}}{\varepsilon}\right) \right\rceil, \qquad \varepsilon > 0.$$

PROOF SKETCH. The main step is to show that $K^{(\rho)}$ is

$$\sqrt{1 - \frac{2\pi\rho_{\min}}{\binom{n}{2}}}\text{-contracting}.$$

We do this as in Lemma 1, showing that for any $x, y \in SO(n)$ with $D(x,y) = r$, there exists a coupling $(X, Y)$ of $(K_x^{(\rho)}, K_y^{(\rho)})$ with

$$\mathbb{E}[D(X,Y)^2] \leq \left(1 - \frac{2\pi\rho_{\min}}{\binom{n}{2}}\right)r^2 + O(r^3).$$

To do this, we first note that $0 \leq 2\pi\rho_{\min} \leq 1$ and write $\rho$ as a *mixture*:

$$\rho = 2\pi\rho_{\min} g + (1 - 2\pi\rho_{\min})h,$$

where $g \equiv 1/2\pi$ is the uniform density and

$$h(\theta) = \frac{\rho(\theta) - \rho_{\min}}{1 - 2\pi\rho_{\min}}$$

is another density. We will set $X = R(i,j,\theta)$, $Y = R(i,j,\theta')$ as in the proof of Lemma 1, choosing $1 \leq i < j \leq n$ uniformly at random. The choices of $\theta, \theta'$ will be made as follows:

1. with probability $2\pi\rho_{\min}$, we pick $\theta$ from the uniform density $g$ and set $\theta' = (\theta - \alpha) \mod 2\pi$ as in the previous proof;
2. with probability $1 - 2\pi\rho_{\min}$, we pick $\theta$ with density $h$ and set $\theta' = \theta$.

Using the notation and reasoning in the previous proof, we immediately see that in the first case $D(X,Y)^2 = \|h\|_{\text{hs}}^2 - \langle h, a_{ij}\rangle_{\text{hs}}^2 + O(r^3)$, whereas in the second case $D(X,Y)^2 = r$. It follows that

$$\mathbb{E}[D(X,Y)^2] = 2\pi\rho_{\min}\left\{\|h\|_{\text{hs}}^2 - \frac{1}{\binom{n}{2}}\sum_{1 \leq i < j \leq n}\langle h, a_{ij}\rangle_{\text{hs}}^2\right\}$$

$$+ (1 - 2\pi\rho_{\min})r^2 + O(r^3)$$

$$= \left(1 - \frac{2\pi\rho_{\min}}{\binom{n}{2}}\right)r^2 + O(r^3). \qquad \square$$

5.2. *A random walk on unitary matrices.* In this section we consider a random walk on unitary matrices. To define it properly, we need a set of definitions analogous to that in Section 4.1, which we briefly state below.

$M(n, \mathbb{C})$ is the set of all complex $n \times n$ matrices. In the present setting $a^*$ is the *conjugate transpose* of $a \in M(n, \mathbb{C})$ and we can define the Hilbert–Schmidt inner product (and corresponding norm) via

$$\langle a, b \rangle_{\mathrm{hs}} \equiv \mathrm{Tr}(ab^*), \qquad a, b \in M(n, \mathbb{C}).$$

With this inner product, $M(n, \mathbb{C})$ is isomorphic to $\mathbb{C}^{n^2}$ with the Euclidean inner product. Call $a \in M(n, \mathbb{C})$ *unitary* if $aa^* = a^*a = \mathrm{id}$, the identity matrix. The set $U(n) \subset M(n, \mathbb{C})$ of all $n \times n$ unitary matrices is a smooth, compact submanifold of $M(n, \mathbb{C})$, which is also a Lie group. The metric $D$ in this case is the Riemanian metric induced on $U(n)$ by the Hilbert–Schmidt inner product on the ambient space $M(n, \mathbb{C})$, which is again invariant by multiplication. Moreover, there exists a multiplication-invariant probability measure on $U(n)$ which we again denote by the Haar measure $\mathcal{H}$.

Let $e_1, \ldots, e_n$ be the canonical basis for $\mathbb{C}^n$. For each $1 \le i < j \le n$ fix a (linear) isometry $I_{ij}: \mathrm{span}\{e_i, e_j\} \to \mathbb{C}^2$. If $u \in U(2)$, we let $u_{ij} \in U(n)$ be the unitary operator that acts as $I_{ij}^{-1} \circ u \circ I_{ij}$ on $\mathrm{span}\{e_i, e_j\}$ and as the identity on $\mathrm{span}\{e_i, e_j\}^\perp$ (that is, $u_{ij}$ acts "like" $u$ on $e_i, e_j$). Our random walk is defined by the kernel $S$ given by

$$L_x(S) \equiv \frac{1}{\binom{n}{2}} \sum_{1 \le i < j \le n} \int_{U(2)} \delta_{R_{ij}x}(S) \, dH(R),$$

where $H$ is the Haar measure on $U(2)$. Thus, $X =_{\mathcal{L}} L_x$ is obtained from $x$ by first choosing $i, j$ uniformly at random, then picking $R \in U(2)$ from the $(2 \times 2)$ Haar measure independently from $i, j$ and then letting $R_{ij}$ act over the two-dimensional subspace $\mathrm{span}\{e_i, e_j\}$.

Our main goal will be to prove an analogue of Theorem 1 in this setting.

THEOREM 5. *Let $\tau_{D,2}(\cdot)$ denote the $L^2$ transportation-cost mixing time for $(M, d) = (U(n), D)$ and $P = L$ as just defined. Then*

$$\tau_{D,2}(\varepsilon) \le \left\lceil n^2 \ln\left(\frac{\pi \sqrt{n}}{\varepsilon}\right) \right\rceil, \qquad \varepsilon > 0.$$

PROOF SKETCH. According to Corollary 1, we need two ingredients: a $\pi\sqrt{n}$ bound on the diameter of $U(n)$ and a "local contraction" estimate for $(L_x, L_y)$ akin to Lemma 1. The diameter bound is easily obtained. Any $u \in U(n)$ has orthogonal eigenvectors with eigenvalues of the form $e^{\sqrt{-1}\theta_i}$ for $\theta_i \in [-\pi, \pi]$, $1 \le i \le n$. For all $t \in [0, 1]$, $u^t \in U(n)$ is a matrix with the same eigenbasis and eigenvalues $e^{\sqrt{-1}t\theta}$, hence, $u^t \in U(n)$ always. The curve $t \mapsto u^t$ ($t \in [0,1]$) has constant speed equal to

$$\sqrt{\sum_{i=1}^n |\theta_i|^2} \leq \pi\sqrt{n}$$

and connects id to $u$; any $x,y$ can be connected by the curve $t \mapsto (yx^*)^t x$, which also has length $\leq \pi\sqrt{n}$, hence, $D(x,y) \leq \pi\sqrt{n}$ for all $x,y \in U(n)$, as desired.

We now provide a local contraction estimate. The key realization is that the tangent space of $U(n)$ at the identity is

$$T = T_{\text{id}}(U(n)) = \{h \in M(n,\mathbb{C}): h = -h^*\}.$$

This means that if $x,y \in U(n)$ and $D(x,y) = r$,

$$\|yx^* - \text{id} - h\|_{\text{hs}} = O(r^2) \qquad \text{for } h = \Pi_T(yx^* - \text{id}),$$

$\Pi_T$ being the orthogonal projector of $M(n,\mathbb{C})$ onto $T$ (as in the previous proof). Moreover, the estimates in Section 4.2 carry over to our current setting.

Suppose $x,y$ as above are given. We choose $1 \leq i < j \leq n$ uniformly at random, $R \in U(2)$ from the Haar measure and will set $R' = Rv$ for some $v = v(i,j,x,y)$ in $U(2)$ to be chosen, so that $R'$ is also Haar distributed on $U(2)$, independently of $i,j,x,y$. This implies that

$$(X,Y) = (R_{ij}x, R'_{ij}y)$$

is a valid coupling of $(L_x, L_y)$. Moreover,

$$D(X,Y) = D(v_{ij}, yx^*).$$

We will now define an orthonormal basis for $M(n,\mathbb{C})$. For $k,\ell \in \{1,\ldots,n\}$, let $u_{k\to\ell}$ be the unique linear operator that maps $e_k$ to $e_\ell$ and $e_t$ to $0$ for all $t \neq k$. The matrices $\{u_{k\to\ell}\}_{1 \leq k\ell \leq n}$ form a orthogonal basis of $M(n,\mathbb{C})$. Since $h^* = -h$, one can check that

$$h = \sum_{k=1}^n \sqrt{-1}h(k,k)u_{k\to k} + \sum_{1 \leq k < \ell \leq n} (h(k,\ell)u_{k\to\ell} - \overline{h(k,\ell)}u_{\ell\to k}),$$

with $h(k,k) \in \mathbb{R}$ and $h(k,\ell) \in \mathbb{C}$. By orthogonality, we have

$$\|h\|_{\text{hs}}^2 = \sum_{k=1}^n h(k,k)^2 + 2\sum_{1 \leq k < \ell \leq n} |h(k,\ell)|^2.$$

We will make a choice of $v$ such that

$$v_{ij} = I_{ij}^{-1} \circ v \circ I_{ij} = e^{h_{ij}} \qquad \text{with}$$

$$h_{ij} \equiv (\sqrt{-1}(h(i,i) + h(j,j)) + h(i,j)u_{i\to j} - \overline{h(i,j)}u_{j\to i}).$$

Indeed, since $h_{ij}^* = -h_{ij}$, $v_{ij} \in U(n)$. Moreover, since $e^{h_{ij}} e_t = e_t$ for $t \neq i, j$, this $e^{h_{ij}}$ acts nontrivially only on span$\{e_i, e_j\}$ and one can easily see that this implies the existence of the desired $v$. Finally, this $v$ only depends on $i, j$ and $x, y$ [through $h = \Pi_T (yx^* - \text{id})$], therefore, it is a valid choice for the coupling construction of $R$ and $R' = Rv$.

One can check that $\|v_{ij} - \text{id}\|_{\text{hs}} = O(r)$, that $\|v - \text{id} - h_{ij}\|_{\text{hs}} = O(r^2)$ and, therefore,

$$\begin{aligned}
D(X, Y)^2 &= D(v_{ij}, yx^*)^2 \\
&= \|v_{ij} - yx^*\|_{\text{hs}}^2 + O(r^3) \\
&= \|h_{ij} - h\|_{\text{hs}}^2 + O(r^3) \\
(\text{expand } h_{ij} - h) &= \|h\|_{\text{hs}}^2 - h(i,i)^2 - h(j,j)^2 - 2h(i,j)^2.
\end{aligned}$$

Averaging over the choices of $u$, $i$ and $j$, we get

$$(18) \quad \mathbb{E}[D(X,Y)^2] = \|h\|_{\text{hs}}^2 - \frac{2}{n}\sum_{i=1}^{n} h(i,i)^2 - \frac{1}{\binom{n}{2}} \sum_{1 \leq i < j \leq n} 2h(i,j)^2 + O(r^3)$$

$$(19) \quad \leq \left(1 - \frac{1}{\binom{n}{2}}\right) r^2 + O(r^3).$$

This implies that the chain $L$ is

$$\sqrt{1 - \frac{1}{\binom{n}{2}}}\text{-locally contracting,}$$

which implies the desired result via Theorem 3. □

**6. Lower bounds for mixing times.** In this section we prove a general mixing time lower bound for random walks induced by group actions. Again, let $(M, d)$ be a metric space.

ASSUMPTION 1. $M$ is compact (hence Polish). There exists a group $G$ acting isometrically on $M$ on the left. That means that there exists a mapping taking $(g, x) \in G \times M$ to $gx \in M$ such that for all $g, h \in G$, $g(hx) = (gh)x$ and for all $g \in G$, $x, y \in M$, $d(gx, gy) = d(x, y)$. We also assume that there is a metric $\tilde{d}$ on $G$ such that $(G, \tilde{d})$ is compact and

$$\forall g, h \in G, \qquad \tilde{d}(g, h) \geq \sup_{x \in M} d(gx, hx).$$

Finally, a Markov transition kernel $P$ on $M$ is defined via a probability measure $\alpha$ on $G$ as follows:

$$\forall x \in M, \qquad P_x \equiv \int_G \delta_{hx}\, d\alpha(h).$$

That is, to sample $X =_{\mathcal{L}} P_x$, one samples $h =_{\mathcal{L}} \alpha$ and sets $X = hx$.

One can check $P$ is indeed a Markov transition kernel; indeed, this follows from the fact that $x \mapsto P_x$ is 1-Lipschitz as a map from $(M,d)$ to $(\Pr(M), W_{d,1})$.

It is well known that compactness of $(M,d)$ and $(G, \tilde{d})$ imply the following (we will use $\sim$ to denote all quantities related to the metric $\tilde{d}$):

- For all $r > 0$ and $H \subset G$, $H$ can be covered by finitely many open balls of radius $r$ in $G$; the minimal number of balls in such a covering is called the $r$-*covering number* of $H$ and denoted by $\tilde{C}_H(r)$.
- For all $r > 0$, there exists a number $N_M(r)$, called the $r$-*packing number* of $M$, which is the largest cardinality of a subset $S \subset M$ with $d(s, s') > r$ for all distinct $s, s' \in S$ (we call such an $S$ *maximally $r$-sparse*).

We can now state our general lower bound result.

THEOREM 6. *Under Assumption 1, suppose that there exists a measure $\mu_* \in \Pr(M)$ and numbers $\tau \in \mathbb{N}$, $\varepsilon > 0$ and $p \geq 1$ such that*

$$\forall x \in M, \qquad W_{d,p}(P_x^\tau, \mu_*) \leq \varepsilon.$$

*Then*

$$\tau \geq \frac{\ln N_M(8\varepsilon) - \ln 2}{\ln \tilde{C}_H(\varepsilon/\tau)},$$

*where $H$ is the support of $\alpha$.*

To understand Theorem 6, it is a good idea to consider the special case $M = G$ is a finite-dimensional Lie group (acting on itself by left-multiplication), $\mu_*$ is a Haar measure on $G$, $P_x^t \to \mu_*$ for all $x \in G$ as $t \to +\infty$ and $\tau = \tau_{d,p}(\varepsilon)$ is the $\varepsilon$-mixing time. Since $G$ is a Lie group, thus a smooth manifold that is locally Euclidean, one would expect that

$$\ln N_G(r) \approx (\text{dimension of } G) \ln(1/r), \qquad 0 < r \ll 1.$$

Similarly, if $H$ has a dimension (in some loosely defined sense), we expect that

$$\ln \tilde{C}_H(r) \approx (\text{dimension of } H) \ln(1/r), \qquad 0 < r \ll 1.$$

Thus, for small enough $\varepsilon$, one would have

$$\tau_{d,p}(\varepsilon) \geq \frac{(\text{dimension of } G)}{(\text{dimension of } H)},$$

at least up to constant factors. The upshot is that a "small" (low-dimensional) set of generators $H$ cannot generate a "large" (high-dimensional) group $G$ in time less than the ratio of the dimensions.

Of course, the reasoning we just presented is not a rigorous proof. In the particular case of Kac's walk, we will need to have bounds on $\tilde{C}_H(\varepsilon)$ and $N_G(\varepsilon)$ that work for a fixed $\varepsilon$, not for $\varepsilon \to 0$.

Let us now prove the general theorem (the bound for Kac's walk is proven subsequently).

PROOF OF THEOREM 6. Let $c = \tilde{C}_H(\varepsilon/\tau)$. By assumption, $H$ can be covered by $c$ open balls of radius $\varepsilon/\tau$ according to $\tilde{d}$, which we represent with $\tilde{B}$:

$$H \subset \tilde{B}(h_1, \varepsilon/\tau) \cup \cdots \cup \tilde{B}(h_c, \varepsilon/\tau).$$

Define the sets $\tilde{S}_1 = \tilde{B}(h_1, \varepsilon/\tau) \cap H$ and $\tilde{S}_i = \tilde{B}(h_i, \varepsilon/\tau) \cap H \setminus \bigcup_{j=1}^{i-1} \tilde{B}(h_j, \varepsilon/\tau)$. These sets form a partition of $H$, hence, the following sum defines a probability measure supported on $\{h_1, \ldots, h_c\}$:

$$\beta \equiv \sum_{i=1}^{c} \alpha(\tilde{S}_i) \delta_{h_i}.$$

In fact, $\beta$ is the image of $\alpha$ under the map $\Psi$ that maps the elements of $\tilde{S}_i$ to $h_i$, for each $i \in \{1, \ldots, c\}$. This map satisfies $\tilde{d}(\Psi(h), h) < \varepsilon/\tau$ because $\tilde{S}_i \subset \tilde{B}(h_i, \varepsilon/\tau)$ by construction. One may check that this implies $W_{\tilde{d},p}(\alpha, \beta) < \varepsilon/\tau$.

Let $Q$ be the Markov transition kernel corresponding to $\beta$ in the same way that $P$ corresponds to $\alpha$; that is,

$$\forall x \in M, \qquad Q_x = \int_{\{h_1,\ldots,h_c\}} \delta_{hx}\, d\beta(h) = \sum_{i=1}^{c} \alpha(\tilde{S}_i) \delta_{h_i x}.$$

For any $x \in M$, if the random pair $(A, B)$ is a coupling of $(\alpha, \beta)$ with $\mathbb{E}[\tilde{d}(A, B)^p]^{1/p} < \varepsilon/\tau$, $(Ax, Bx)$ is a coupling of $(P_x, Q_x)$ with

$$W_{d,p}(P_x, Q_x) \leq \mathbb{E}[d(Ax, Bx)^p]^{1/p} \leq \mathbb{E}[\tilde{d}(A, B)^p]^{1/p} < \varepsilon/\tau.$$

Hence,

$$\forall x \in M, \qquad W_{d,p}(P_x, Q_x) < \varepsilon/\tau.$$

A simple calculation implies

$$\forall x \in M, \qquad W_{d,p}(Q_x^\tau, P_x^\tau) < \varepsilon.$$

For any $x \in M$, the definition of $\tau$ implies $W_{d,p}(\mu_*, P_x^\tau) \leq \varepsilon$, so that

$$W_{d,p}(Q_x^\tau, \mu_*) \leq W_{d,p}(P_x^\tau, \mu_*) + W_{d,p}(P_x^\tau, Q_x^\tau) < 2\varepsilon.$$

Thus, the $p$-optimal coupling $(X_x, Y)$ of $(Q_x^\tau, \mu_*)$ achieves

$$\mathbb{E}[d(X_x, Y)^p]^{1/p} < 2\varepsilon.$$

Now let $S \subset M$ be a maximal $8\varepsilon$-sparse subset of $M$. Notice that the cardinality of $S$ is $\ell \equiv N_M(8\varepsilon)$, by definition of the latter quantity. One may define random variables $\{X_x\}_{x \in S}, Y$ on the same probability space such that, for each $x \in S$, $(X_x, Y)$ is a coupling of $(Q_x^\tau, \mu_*)$ achieving the above bound (this follows, e.g., from the Gluing lemma in the Introduction of Villani's book [15]). Hence, if

$$I_x = \chi_{\{d(X_x, Y) \geq 4\varepsilon\}} \qquad (x \in S),$$

Markov's inequality implies that

$$\mathbb{P}(I_x = 1) < 1/2$$

and

$$\mathbb{E}\left[\sum_{x \in S} I_x\right] < \frac{\ell}{2}.$$

It follows that there *exists* a realization of $\{X_x\}_{x \in S}, Y$ such that $d(X_x, Y) < 4\varepsilon$ for all $x$ in a subset $S' \subset S$ of cardinality $\geq \ell/2$.

We *fix* such a realization. For each $x \in S$, the support of the measure $Q_x$ is contained in the finite set $\{h_1 x, \ldots, h_c x\}$. A simple inductive argument shows that $X_x = v_x x$ for some

$$v_x \in V_\tau \equiv \{h_{i_1} h_{i_2} \cdots h_{i_\tau} : i_1, i_2, \ldots, i_\tau \in \{1, 2, \ldots, c\}\}.$$

Now notice that, on the one hand, for all $x, x' \in S'$,

$$d(X_x, X_{x'}) \leq d(X_x, Y) + d(X_{x'}, Y) < 4\varepsilon + 4\varepsilon = 8\varepsilon.$$

On the other hand, for distinct $x, x' \in S$, if $v_x = v_{x'}$, then $d(X_x, X_{x'}) = d(v_x x, v_x x') = d(x, x') \geq 8\varepsilon$ since $S$ is $8\varepsilon$-sparse and $d$ is invariant by left multiplication. We deduce that

$$\forall x, x' \in S', \qquad x \neq x' \quad \Rightarrow \quad v_x \neq v_{x'}.$$

This implies that

$$\ell/2 \leq \text{cardinality of } S' \leq \text{cardinality of } V_\tau$$

and the latter quantity is clearly upper bounded by $c^\tau$. We deduce that

$$\ell/2 \leq c^\tau \quad \Rightarrow \quad \tau \geq \frac{\ln \ell - \ln 2}{\ln c}.$$

The proof is finished once we recall that $\ell = N_M(8\varepsilon)$ and $c = \tilde{C}_H(\varepsilon/\tau)$. $\square$

6.1. *The lower bound for Kac's random walk (Theorem 2).* We now show that Theorem 2 follows from the general lower bound.

PROOF OF THEOREM 2. We will freely use the notation introduced in Section 4. In particular, we take $M = G = SO(n)$, $d = \mathrm{hs}$, $P = K$, $\mu_* = \mathcal{H}$, $\tau = \tau_{d,p}(\varepsilon)$ and

$$H = \bigcup_{1 \leq i < j \leq n} \{R(i,j,\theta) : \theta \in [0, 2\pi]\} \tag{20}$$

and

$$\alpha = \frac{1}{\binom{n}{2}} \sum_{1 \leq i < j \leq n} \int_0^{2\pi} \delta_{R(i,j,\theta)} \, d\theta.$$

Notice that $d$ is right-invariant and that we may take $\tilde{d} = d$ in this case.

We now upper bound $\tilde{C}_H(r)$. Equation (20) shows that $H$ is the union of $\binom{n}{2}$ sets. Each of those is an isometric image of the unit circle in the intrinsic metric $D$ of $SO(n)$. Since hs is dominated by $D$, we have

$$\forall 0 < r \leq 2\pi, \qquad \tilde{C}_H(r) \leq 2\pi \binom{n}{2} r^{-1} \leq \pi n^2 r^{-1}. \tag{21}$$

We must also lower bound $N_M(r)$. A maximal $r$-packing $S \subset SO(n)$ has to satisfy $\min_{s \in S} d(x,s) < r$ for all $x \in SO(n)$ (an $x$ violating the bound could be added to $S$, which violates maximality). This implies that

$$\begin{aligned} SO(n) = \bigcup_{s \in S} B(s,r) &\Rightarrow \sum_{s \in S} \mathcal{H}(B(s,r)) \geq 1 \\ &\Rightarrow_{\text{(Inv)}} \quad N_M(r) = |S| \geq \frac{1}{\mathcal{H}(B(\mathrm{id},r))}. \end{aligned} \tag{22}$$

[Implication (Inv) uses the invariance of $\mathcal{H}$, which implies that all balls of radius $r$ have the same measure.]

We now make the following claim:

CLAIM 1. *There exist constants $\phi, \psi > 0$ such that, for all $n \geq 10$ and $0 < r < \sqrt{n}/10$,*

$$\mathcal{H}(B(\mathrm{id}, r)) \leq \left(\frac{e^\phi r}{\sqrt{n}}\right)^{\psi n^2}.$$

The restrictions on $n, r$ are by no means sharp, but they give us some room to spare in what follows.

Before proving the claim, we show how the theorem follows from it. Given $n \geq 10$, assume $\tau = \tau_{d,p}(\varepsilon) \leq n^3/\pi$. We see that, for $0 < \varepsilon < 1$,

$$\tau \ln \tilde{C}_H(\varepsilon/\tau) \leq (\ln \tau + 2\ln n + \ln \pi + \ln(1/\varepsilon))\tau$$
$$\leq (5\ln n + \ln(1/\varepsilon))\tau \qquad (\text{use } \tau \leq n^3/\pi);$$
$$\ln N_M(8\varepsilon) \geq -\ln \mathcal{H}(\text{id}, 8\varepsilon) \qquad [\text{via (21)}]$$
$$\geq \psi n^2 \left(\frac{\ln n}{2} + \ln(1/8\varepsilon) - \phi\right).$$

This implies that, for $\varepsilon \equiv e^{-\phi}/8 < 1$, one can use the bound in Theorem 6 to see that

$$\tau \geq \frac{\psi/2 n^2 \ln n}{5\ln n + \phi + \ln 8} \geq cn^2$$

for some $c > 0$ not depending on $n$. Of course, if $\tau > n^3/\pi$, $\tau > cn^2$ for a (possibly smaller) $c > 0$, so the inequality presented above actually implies the theorem for $n \geq 10$. Since there is only a finite set of remaining values of $n$, one may finish the proof by picking a smaller $c$, if necessary.

It remains to prove the claim. We will do so via probabilistic reasoning, using some rough upper estimates and known results for spheres in an arbitrary dimension. As a preliminary, consider $x \in SO(n)$ and let $x_i \in \mathbb{R}^n$ denote its $i$th column. One has

$$\|x - \text{id}\|_{\text{hs}}^2 = \sum_{i=1}^n |x_i - e_i|^2.$$

The columns of $x$ are orthonormal, hence, $|x_i| = |e_i| = 1$ and

$$\|x - \text{id}\|_{\text{hs}}^2 = 2\sum_{i=1}^n (1 - x_i.e_i).$$

Hence,

$$\|x - \text{id}\|_{\text{hs}} < r \quad \Rightarrow \quad \frac{1}{n}\sum_{i=1}^n (1 - x_i.e_i) < \frac{r^2}{2n}.$$

One can now use Markov's inequality to deduce that

(23)
$$\|x - \text{id}\|_{\text{hs}} < r \quad \Rightarrow \quad \exists I \subset \{1, \ldots, n\} \qquad \text{with } |I| = \lceil n/2 \rceil \text{ such that}$$
$$\forall i \in I, \qquad x_i.e_i > 1 - \frac{r^2}{n}.$$

Thus, if $X =_{\mathcal{L}} \mathcal{H}$ is a random variable, defined on some probability space $(\Omega, \mathcal{F}, \mathbb{P})$ and with values in $SO(n)$,

$$\mathcal{H}(B(\text{id}, r)) = \mathbb{P}(\|X - \text{id}\| < r)$$

$$\leq \sum_{I \subset \{1,\ldots,n\}\,:\,|I|=\lceil n/2 \rceil} \mathbb{P}\left(\forall i \in I, X_i.e_i > 1 - \frac{r^2}{n}\right)$$

(24) $\qquad (\mathcal{L}_{(X_i : i \in I)}$ is the same for all $I$)

$$= \binom{n}{\lceil n/2 \rceil} \mathbb{P}\left(\forall 1 \leq i \leq \lceil n/2 \rceil, X_i.e_i > 1 - \frac{r^2}{n}\right)$$

$$\leq 2^n \mathbb{P}\left(\forall 1 \leq i \leq \lceil n/2 \rceil, X_i.e_i > 1 - \frac{r^2}{n}\right).$$

Now consider the orthogonal projection maps:

$$\Pi_i = \Pi_i(X): \qquad z \in \mathbb{R}^n \mapsto z - \sum_{k=1}^{i-1}(X_k, z)X_k,$$

with $\Pi_1$ is the identity operator. Clearly, $0 < \Pi_i e_i < 1$ for all $1 \leq i \leq \lceil n/2 \rceil$ with probability 1. $X_i$ belongs to the range of $\Pi_i$, a self-adjoint operator. Hence, outside of a null set,

$$X_i.e_i = \Pi_i X_i.e_i = X_i.\Pi_i e_i > 1 - \frac{r^2}{n} \quad \Rightarrow \quad X_i.\frac{\Pi_i e_i}{|\Pi_i e_i|} > 1 - \frac{r^2}{n}.$$

This implies the bound

$$\mathbb{P}\left(\forall 1 \leq i \leq \lceil n/2 \rceil, X_i.e_i > 1 - \frac{r^2}{n}\right) \leq \mathbb{P}\left(\bigcap_{i=1}^{\lceil n/2 \rceil} E_i\right)$$

$$\left(\text{with } E_i \equiv \left\{X_i.\frac{\Pi_i e_i}{|\Pi_i e_i|} > 1 - \frac{r^2}{n}\right\}\right).$$

Let $\mathcal{F}_0 = \{\varnothing, \Omega\}$ be the trivial $\sigma$-field on $\Omega$ and, for $1 \leq j \leq n$, $\mathcal{F}_j$ be the $\sigma$-field generated by $X_1, \ldots, X_j$. These $\sigma$-algebras form an increasing sequence. We omit the proof of the following three facts, valid for each $1 \leq i \leq \lceil n/2 \rceil$:

1. $E_i$ is $\mathcal{F}_i$-measurable;
2. $\Pi_i e_i / |\Pi_i e_i|$ $\mathcal{F}_{i-1}$-measurable;
3. conditioned on $\mathcal{F}_{i-1}$, $X_i$ is uniform on the $(n-i)$-dimensional unit sphere of the subspace of $\mathbb{R}^n$ corresponding to the range of $\Pi_i$, and $\Pi_i e_i / |\Pi_i e_i|$ is a point on that same sphere.

Let $\mathcal{V}_{n-i}$ be the (normalized) uniform measure on $S^{n-i}$. The above considerations, together with the rotational invariance of $\mathcal{V}_{n-i}$, imply that

$$\forall 1 \leq i \leq \lceil n/2 \rceil, \qquad \mathbb{P}(E_i \mid \mathcal{F}_{i-1}) = \mathcal{V}_{n-i}(C_{n-i}(1 - r^2/n)) \qquad \text{a.s.,}$$

where, for a given $\tau \in \mathbb{R}$, $C_{n-i}(\tau)$ is the spherical cap

$$C_{n-i}(\tau) \equiv \{v \in S^{n-i} : v.e_1 > \tau\}.$$

A simple inductive argument with conditional expectations then shows that

$$\mathbb{P}\left(\bigcap_{i=1}^{\lceil n/2 \rceil} E_i\right) = \mathbb{E}[\mathbb{P}(E_{\lceil n/2 \rceil} \mid \mathcal{F}_{\lceil n/2 \rceil - 1}) \chi_{\bigcap_{i=1}^{\lceil n/2 \rceil - 1} E_i}]$$

$$= \mathcal{V}_{n-\lceil n/2 \rceil}(C_{n-\lceil n/2 \rceil}(1 - r^2/n)) \mathbb{P}\left(\bigcap_{i=1}^{\lceil n/2 \rceil - 1} E_i\right)$$

$$= (\cdots)$$

$$= \prod_{i=1}^{\lceil n/2 \rceil} \mathcal{V}_{n-i}(C_{n-i}(1 - r^2/n)).$$

We now apply known bounds on the volume of spherical caps [2], Lemma 2.1:

(25)
$$\forall m \in \mathbb{N} \setminus \{0, 1\}, \forall \tau \in [2/\sqrt{m}, 1],$$
$$\frac{(1-\tau^2)^{(m-1)/2}}{6\tau\sqrt{m}} \leq \mathcal{V}_m(C_m(\tau)) \leq \frac{(1-\tau^2)^{(m-1)/2}}{2\tau\sqrt{m}}.$$

We need the upper bound for $n - \lceil n/2 \rceil \leq m \leq n - 1$ and

$$\tau = 1 - \frac{r^2}{n}, \quad \text{which} \in \left[\sqrt{\frac{2}{m}}, 1\right] \quad \text{for } n \geq 10, r \leq \sqrt{n}/10.$$

Moreover, we know that in this case $2\tau^2 > 2 - 4r^2/n > 1$, so

$$\mathbb{P}\left(\bigcap_{i=1}^{\lceil n/2 \rceil} E_i\right) \leq \prod_{i=1}^{\lceil n/2 \rceil} \mathcal{V}_{n-i}(C_{n-i}(1 - r^2/n)) \leq \prod_{i=1}^{\lceil n/2 \rceil} \frac{(2r^2)^{(n-i)/2}}{\sqrt{n-i}} \leq \left(\frac{e^{\phi_0} r}{\sqrt{n}}\right)^{\psi n^2},$$

for some constants $\phi_0, \psi > 0$ not depending on $n \geq 10$ or $r \leq \sqrt{n}/10$. Using (24), we deduce that

$$\mathcal{H}(B(\mathrm{id}, r)) \leq 2^n \mathbb{P}\left(\bigcap_{i=1}^{\lceil n/2 \rceil} E_i\right) \leq \left(\frac{e^{\phi} r}{\sqrt{n}}\right)^{\psi n^2},$$

with $\phi > 0$ another constant. The claim and the theorem are finally proven. □

## 7. Final remarks.

- The most obvious problem left open in the present paper is a sharp characterization of the mixing time of Kac's walk. We conjecture that our upper bound is tight for all $\varepsilon \in (0, \varepsilon_0)$; that is, that there exist constants $c, \varepsilon_0 > 0$ such that, for all $n \geq 3$ and $\varepsilon \in (0, \varepsilon_0)$,

$$\tau_{\mathrm{hs},1}(\varepsilon) \geq cn^2 \ln\left(\frac{cn}{\varepsilon}\right).$$

Notice that the restriction to $n \geq 3$ is necessary, as for $n = 2$ the walk mixes perfectly in one single step.

The quantity $n^2 \ln n$ in the conjectured lower bound immediately suggests a "coupon-collector phenomenon." For instance, one is tempted to guess that the walk cannot mix before 2-dimensional rotations have been applied to all possible pairs $e_i, e_j$ of canonical basis vectors. The difficulty with this idea is that two rows of $X(t)$ may "interact" without ever being changed in the same step of the walk.

- The simple lower bound method in Section 6 cannot go farther than $\Omega(n^2)$, even if $\varepsilon \to 0$ with $n$. It would be interesting to derive better lower bounds at this level of generality.
- Going back to the application of Ailon and Chazelle [1], $O(n^2 \ln n)$ mixing is still too large for $n$ big, which is precisely when dimensionality reduction is the most useful. However, that application only requires that certain projections behave as they should, which is a less stringent requirement than approximating the Haar measure. It is thus natural to ask whether better bounds might be available for that specific application. More precisely, let $Y_k(t) = \Pi_k X(t)^\dagger$, where $X(t)$ is a realization of Kac's walk and $\Pi_k$ is the projection onto the first $k$ canonical basis vectors. Clearly, $\{Y_k(t)\}_{t=0}^{+\infty}$ corresponds to a Markov chain on the *Stiefel manifold*:

$$V_k(\mathbb{R}^n) \equiv \{(v_1, \ldots, v_k) \in (\mathbb{R}^n)^k : \forall 1 \leq i, j \leq k, v_i.v_j = \delta_{ij}\}.$$

One can adapt the proof of Theorem 2 to show that this walk cannot mix in less than $\Omega(nk)$ time.

We conjecture that $Y_k(t)$ mixes in $\Theta(nk \ln n)$ steps. Recall that for dimension reduction we need $k = O(\ln |S|)$. Our conjecture would imply great time savings for $n \gg \ln |S|$.

- Theorem 3 on local-to-global coupling can be used to reprove some known results. Consider, for instance, a Riemannian manifold $M$ with dimension $n$, distance $d$ and Ricci curvature lower bounded by $K \in \mathbb{R}$. Let $P = P^{(\varepsilon)}$ correspond to the ball walk on $M$ where a step from $x$ consists of choosing $X$ uniformly from the ball $B(x, \varepsilon)$. Using a simple, "strictly local" variant of [17], Lemma 2, and our Theorem 3, one can very easily show that $\mu \mapsto \mu P^{(\varepsilon)}$ is $(1 - K\varepsilon^2/2(n+2) + o(\varepsilon^2))$-Lipschitz (thus contracting when $K > 0$ and $\varepsilon$ is small enough). By "strictly local," we mean that we do *not* need to have control $W_{d,1}(P_x, P_y)$ uniformly over all pairs of nearby points in the manifold: we just need that for each *fixed* $x \in M$, as $y \to x$,

$$W_{d,p}(P_y^{(\varepsilon)}, P_x^{(\varepsilon)}) \leq (\xi + o(1))d(x, y)$$

for the appropriate $\xi > 0$.

We expect that checking the local Lipschitz condition in other applications will oftentimes be much simpler than proving a global contraction estimate.

**Acknowledgment.** We thank Yann Ollivier for interesting discussions on Ricci curvature and Markov chains. We are also grateful to Christian Rau, who found several typos and minor errors in a previous version of the present work.

Instituto de Matemática Pura
e Aplicada (IMPA)
Estrada Dona Castorina, 110
Rio de Janeiro, RJ 22460-320
Brazil
E-mail: rimfo@impa.br